\documentclass[10pt]{amsart}
\usepackage{amsmath}
\usepackage{amssymb}
\usepackage{enumerate}
\usepackage{amsbsy}
\usepackage{amsfonts}
\usepackage{color}

\newtheorem{thm}{Theorem}[section]
\newtheorem{lem}[thm]{Lemma}

\numberwithin{equation}{section}

\begin{document}

	\title[Yang-Mills System]{A Note on the smoothness of flow maps for the Yang-Mills system in the Lorenz gauge}

    \author[S. Hong]{Seokchang Hong}
    \address{Department of Mathematical Sciences, Seoul National University, Seoul 08826, Republic of Korea}
    \email{seokchangh11@snu.ac.kr}

	
	\begin{abstract}
		We study the failure of the smoothness of flow maps for the $(1+3)$ dimensional Yang-Mills system in the Lorenz gauge by Knapp type counterexamples. This shows a gap between the scaling critical regularity exponents and the best attainable regularity via Picard's iteration in the Yang-Mills system under the Lorenz gauge condition. 
	\end{abstract}

		\maketitle

\section{Introduction}
In this note, we prove the failure of the $C^2$ smoothness of flow maps for the Yang-Mills system in the Lorenz gauge, which is given by the following nonlinear wave system \cite{hong}: 
\begin{align}\label{ym-system}
\begin{aligned}
\Box A_\nu = & -2[A^\mu,\partial_\mu A_\nu]+[A^\mu,\partial_\nu A_\mu]-[A^\mu,[A_\mu,A_\nu]], \\
\Box F_{\mu \nu}=&-2\left[A^{\lambda}, \partial_{\lambda} F_{\mu \nu}\right]+2\left[\partial_{\nu} A^{\lambda}, \partial_{\lambda} A_{\mu}\right]-2\left[\partial_{\mu} A^{\lambda}, \partial_{\lambda} A_{\nu}\right] \\
&+2\left[\partial^{\lambda} A_{\mu}, \partial_{\lambda} A_{\nu}\right]+2\left[\partial_{\mu} A^{\lambda}, \partial_{\nu} A_{\lambda}\right]-\left[A^{\lambda},\left[A_{\lambda}, F_{\mu \nu}\right]\right] \\
&+2\left[F_{\lambda \mu},\left[A^{\lambda}, A_{\nu}\right]\right]-2\left[F_{\lambda \nu},\left[A^{\lambda}, A_{\mu}\right]\right]-2\left[\left[A^{\lambda}, A_{\mu}\right],\left[A_{\lambda}, A_{\nu}\right]\right],
\end{aligned}
\end{align}
Let $G$ be a compact Lie group and $\mathfrak g$ its Lie algebra.
The gauge potential $A$ is $\mathfrak g$-valued $1$-form and the curvature $2$-form $F$ is given by
$$
F_{\mu\nu} = \partial_\mu A_\nu-\partial_\nu A_\mu + [A_\mu,A_\nu].
$$
Here $[\cdot,\cdot]$ is the matrix commutator. From now on, we assume $G=SU(n,\mathbb C)$ (the group of unitary matirices of determinant one) and $\mathfrak g=\mathfrak{su}(n,\mathbb C)$ (the algebra of trace-free skew-Hermitian matrices).  
We denote the initial data for $(A,F)$ at $t=0$ by
\begin{align}\label{initial-data}
A(0) = a, \quad \partial_tA(0) = \dot{a},\quad F(0) = f,\quad \partial_tF(0) = \dot{f}.
\end{align}
We note that the initial data $f$ for $F$ is completely determined by $(a,\dot{a})$:
\begin{align}\label{initial-f}
\left\{\begin{array}{l}f_{i j}=\partial_{i} a_{j}-\partial_{j} a_{i}+\left[a_{i}, a_{j}\right], \\ 
f_{0 i}=\dot{a}_{i}-\partial_{i} a_{0}+\left[a_{0}, a_{i}\right], \\ 
\dot{f}_{i j}=\partial_{i} \dot{a}_{j}-\partial_{j} \dot{a}_{i}+\left[\dot{a}_{i}, a_{j}\right]+\left[a_{i}, \dot{a}_{j}\right], \\ 
\dot{f}_{0 i}=\partial^{j} f_{j i}+\left[a^{\alpha}, f_{\alpha i}\right].\end{array}\right.
\end{align}
and we also have
\begin{align}\label{initial-const}
\left\{\begin{array}{l}\dot{a}_0 = \partial^ia_i,\\
\partial^if_{0i} = [a^i,\dot{a}_i].\end{array}\right.
\end{align}
In \cite{hong}, the present author studied the local well-posedness of \eqref{ym-system}. 
\begin{thm}[Theorem 1.1 of \cite{hong}]\label{lwp}
Let $\sigma\ge1$. Suppose that given initial data $(a,\dot{a},f,\dot{f})$ which satisfy \eqref{initial-data} and \eqref{initial-const} is in the function space:
$$
\langle\Omega\rangle^{-\sigma} B^\frac12_{2,1}\times\langle\Omega\rangle^{-\sigma} B^{-\frac12}_{2,1}\times\langle\Omega\rangle^{-\sigma} B^{-\frac12}_{2,1}\times\langle\Omega\rangle^{-\sigma} B^{-\frac32}_{2,1},
$$ 
where 
$\langle\Omega\rangle^\sigma$ is the $\sigma$-angular derivative.
Then there exists local existence time $T>0$ depending on the initial data norm
such that there exist solutions 
\begin{align}
\left\{\begin{array}{l}
A\in C([-T,T];\langle\Omega\rangle^{-\sigma}B^\frac12_{2,1})\cap C^1((-T,T);\langle\Omega\rangle^{-\sigma}B^{-\frac12}_{2,1}),\\
 F\in C([-T,T];\langle\Omega\rangle^{-\sigma}B^{-\frac12}_{2,1})\cap C^1((-T,T);\langle\Omega\rangle^{-\sigma}B^{-\frac32}_{2,1}),\end{array}\right.
\end{align}
of \eqref{ym-system}, respectively. 
\end{thm}
In \cite{hong}, we revealed the null structure in $[ A^\mu,\partial_\mu A_\nu]$ and $[A^\mu,\partial_\nu A_\mu]$. The first bilinear term turns out to be the $Q$-type null form and hence the treatment is not difficult. On the other hand, the second bilinear term is problematic. 
\begin{lem}[Lemma 3.5 of \cite{hong}]
In the Lorenz gauge $\partial^\mu A_\mu=0$, we have the following identity:
\begin{align*}
[A^\lambda,\partial_\mu A_\lambda] = \sum_{j=1}^3\Gamma_\mu^j(A,\partial A,F,\partial F),
\end{align*}
where
\begin{align*}
\Gamma^1_\mu(A,\partial A,F,\partial F) &= -[A_0,\partial_\mu A_0]+[D^{-1}R_j(\partial_tA_0),D^{-1}R^j\partial_t(\partial_\mu A_0)] ,\\
\Gamma^2_\mu(A,\partial A,F,\partial F) &= -\frac12\epsilon^{ijk}\epsilon_{klm}(Q_{ij}[D^{-1}R^nA_n,D^{-1}R^l\partial_\mu A_m] \\
&\qquad\qquad\qquad + Q_{ij}[D^{-1}R^n\partial_\mu A^n,D^{-1}R^lA_m]),\\
\Gamma^3_\mu(A,\partial A,F,\partial F) &= [D^{-2}\nabla\times\mathbf F,D^{-2}\nabla\times\partial_\mu\mathbf F] \\
&-[D^{-2}\nabla\times\mathbf F,D^{-2}\partial_\mu\nabla\times(\mathbf A\times\mathbf A)] \\
&-[D^{-2}\nabla\times(\mathbf A\times\mathbf A),D^{-2}\nabla\times\partial_\mu\mathbf F]\\
&+[D^{-2}\nabla\times(\mathbf A\times\mathbf A),D^{-2}\partial_\mu\nabla\times(\mathbf A\times\mathbf A)].
\end{align*}
Here $\mathbf F=(F_{23},F_{31},F_{12})$. 
\end{lem}
The $\Gamma^1,\Gamma^2$ are essentially $Q$-type null forms. From now on, we focus on the $\Gamma^3$. In \cite{hong}, we observe that for initial data $(a,f)\in H^s\times H^r$, the regularity we can attain is at most $r=-\frac14$. More explicitly, if we only use bilinear estimates, then we have the relation between $s$ and $r$: $s-1\le2r$, which is significantly far from $s-r=1$. (See \cite[Remark 6.1]{hong} for details.) From this observation, we expect the Picard's iteration to fail if $s-1>2r$. Indeed, we have the following theorem:
\begin{thm}\label{illp}
For $s,r\in\mathbb R$, suppose that $s-1>2r$. Let $T>0$ be given. Then the flow map of $(a,f)\mapsto (A(t),F(t))$ from $H^s\times H^{r}$ to $C([-T,T];H^s\times H^r)$ cannot be $C^2$ at the origin. 
\end{thm}
This shows the threshold for the contraction mapping principle and reinforce our previous study \cite{hong}. To prove Theorem \ref{illp}, we give explicit counterexample which is Knapp type and follow the argument by \cite{herrlenz,molsatz}. We also refer the readers to \cite{herr,okamoto} for study on regularity of flow maps for other dispersive equations.
In the rest of the paper, we shall introduce the decomposition of d'Alembertian and then give the proof of Theorem \ref{illp}.

\noindent\textbf{Notations.} 
As usual different positive constants dependent on only dimension is denoted by the same letter $C$, if not specified. $A \lesssim B$ and $A \gtrsim B$ means that $A \le CB$ and
$A \ge C^{-1}B$, respectively for some $C>0$. $A \approx B$ means that $A \lesssim B$ and $A \gtrsim B$.

The spatial Fourier transform on $\mathbb{R}^{3}$ is defined by
$$
\widehat{f}(\xi)=\int_{\mathbb{R}^{3}}e^{-ix\cdot\xi}f(x)\,dx,
$$
where $\xi\in\mathbb{R}^{3}$. Also we denote $\mathcal{F}_x(u)=\widehat{u}$. 

We put $D := |\nabla|$ whose symbol is $|\xi|$.
\section{Preliminaries}

\subsection{Decomposition of d'Alembertian}

Given nonlinear wave equation $\Box u = G(u,\partial u)$ with initial data $(u(0),\partial_tu(0))$, we consider the following transform:
$
(u,\partial_tu) \rightarrow (u_+,u_-),
$
where
$$
u_\pm = \frac12\left(u\pm\frac{1}{iD}\partial_tu\right).
$$
Then given nonlinear wave equation is rewritten as
$$
(i\partial_t\pm D)u_\pm = \mp\frac{1}{2D}\widetilde{G}(u_+,u_-),
$$
where $\widetilde{G}(u_+,u_-)=G(u,\partial u)$. The Duhamel integral is given by
\begin{align}\label{duhamel}
u_\pm(t) = u^{\rm hom}_\pm \pm \int_0^te^{\mp i(t-t')D}\frac{1}{2iD}\widetilde{G}(u_+,u_-)(t')\,dt',
\end{align}
where
\begin{align}\label{hom-eq}
u^{\rm hom}_\pm (t) = \frac12e^{\mp itD}\left(u(0)\pm\frac{1}{iD}\partial_tu(0)\right).
\end{align}
\section{The failure of smoothness: Proof of Theorem \ref{illp}}
We prove Theorem \ref{illp} by contradiction. We assume $A_0(t=0)=0$, $A_j(t=0)=\delta a_j$, and $\partial_tA_j(t=0)=0$, where $0<\delta\ll1$. We shall denote the local solutions to \eqref{ym-system} by $(A_\mu(\delta,t),F_{\mu\nu}(\delta,t))$.

\subsection{Set up}
By \eqref{hom-eq} and \eqref{initial-f}, we have
\begin{align}
 &  \sum_\pm F^{\rm hom}_{ij,\pm} = \frac12\sum_\pm e^{\mp itD}\left(f_{ij}\pm\frac{1}{iD}\dot{f}_{ij} \right) = \frac12\sum_\pm e^{\mp itD}\left( \delta\partial_ia_j-\delta\partial_ja_i+\delta^2[a_i,a_j]\right), 
\end{align}
Now we write
\begin{align*}
\partial_\delta^2\Gamma_\mu^3(A,\partial A,F,\partial F)|_{\delta=0} & = 2[D^{-2}\nabla\times\partial_\delta\mathbf F|_{\delta=0},D^{-2}\nabla\times\partial_\mu\partial_\delta\mathbf F|_{\delta=0}], 
\end{align*}
since the derivative of cubic and quartic term in $\Gamma^3$ vanish at $\delta=0$. Thus we recall \eqref{duhamel} to get (we ignore the first bilinear term $[A^\mu,\partial_\mu A_\nu]$.)
\begin{align}
\partial_\delta^2A_\nu|_{\delta=0} & = 2\sum_{\pm_1,\pm_2,\pm_3}\int_0^te^{\mp_1i(t-t')D}(2iD)^{-1}[D^{-2}\nabla\times\partial_\delta\mathbf F_{\pm_2}(t')|_{\delta=0},D^{-2}\nabla\times\partial_\nu\partial_\delta\mathbf F_{\pm_3}(t')|_{\delta=0}]\,dt'.
\end{align}
Hence if the flow is $C^2$ in $H^s\times H^r$ within a local existence time interval $[0,T]$, we have the following estimates:
\begin{align}\label{fail-est}
\sup_{0\le t\le T}\|\partial_\delta^2A_\nu(t)|_{\delta=0}\|_{H^s} &\lesssim \|\mathbf F\|_{H^r}^2.
\end{align}
We note that
\begin{align*}
\partial_\delta F_{ij}(t)|_{\delta=0} & = \frac12\sum_{\pm} e^{\mp itD}(\partial_ia_j-\partial_ja_i). 
\end{align*}
We further assume that $a_3=0$ and $\partial_3a_2=0$ for simplicity. That is, $\widehat{a_2}(\xi)=0$ for all $\xi_3\neq0$, where $\xi=(\xi_1,\xi_2,\xi_3)$. Then
$$
\nabla\times\partial_\delta\mathbf F(t)|_{\delta=0} = \frac12\sum_\pm e^{\mp itD}(\partial_1\partial_2a_2-\partial_2^2a_1,-\partial_1^2a_2+\partial_1\partial_2a_1,0)+\frac12\sum_\pm e^{\mp itD}(-\partial_3^2a_1,0,\partial_1\partial_3a_1).
$$ 
We also assume $a_1=a_1^1T^1$, $a_2=a_2^2T^2$, where $T^a$ are infinitesimal generators of $\mathfrak{su}(n,\mathbb C)$, $a=1,2,\cdots,n^2-1$. Fron now on, we simply put $a_1^1=a_1$ and $a_2^2=a_2$. We let $\lambda\gg1$. Define $W_\lambda = \{ \xi=(\xi_1,\xi_2,\xi_3) : |\xi_1-\lambda|\le10^{-6}\lambda, 10^{-9}\lambda^\frac12\le\xi_2,\xi_3\le10^{-6}\lambda^\frac12 \}$ and $cW_\lambda=\{ c\xi : \xi\in W_\lambda \}$ for $c\in\mathbb R$. We also define $W_\lambda' = \{ \xi = (\xi_1,\xi_2,0) : |\xi_1-\lambda|\le10^{-6}\lambda, |\xi_2|\le 10^{-6}\lambda^\frac12 \}$. Now we take
\begin{align}\label{support-condition}
\widehat{a_1}(\xi) = \chi_{2W_\lambda}(\xi),\quad \widehat{a_2}(\xi) = \chi_{-W_\lambda'}(\xi).
\end{align}
We note that the integration of $\widehat{a_2}$ is $0$, whereas the integration of $\widehat{a_1a_2}$ is nonzero. 

\subsection{Failure of \eqref{fail-est} for $s-1>2r$}
Since $\nabla\times\mathbf F^{\rm hom}$ consists of so many integrations, at first glance it seems tremendous. However, by the support condition \eqref{support-condition} it is no harm to ignore the most of integrations. In fact, for $\xi\in W_\lambda$, we have $\xi_1\approx\lambda$ whereas $\xi_2,\xi_2\approx\lambda^\frac12$. Then we can write $\partial_\delta^2 A_\nu=\Lambda^1_\nu+\Lambda^2_\nu+\cdots$, where the leading terms $\Lambda^1_\nu,\Lambda^2_\nu$ are given by 
\begin{align}
\Lambda^1_\nu & = \frac12\sum_{\pm_1,\pm_2,\pm_3}\int_0^t e^{\mp_1i(t-t')D}(2iD)^{-1}(D^{-2}e^{\mp_2it'D}\partial_1^2a_2D^{-2}\partial_\nu e^{\mp_3it'D}\partial_1\partial_2a_1 \\ 
&\qquad\qquad\qquad\qquad\qquad\qquad\qquad -D^{-2}e^{\mp_2it'D}\partial_1\partial_2a_1D^{-2}\partial_\nu e^{\mp_3it'D}\partial_1^2a_2)\,dt'[T^1,T^2], \nonumber \\
\Lambda^2_\nu & = \frac12\sum_{\pm_1,\pm_2,\pm_3}\int_0^t e^{\mp_1i(t-t')D}(2iD)^{-1}(D^{-2}e^{\mp_2it'D}\partial_1^2a_2D^{-2}\partial_\nu e^{\mp_3it'D}\partial_1\partial_3a_1 \\
&\qquad\qquad\qquad\qquad\qquad\qquad\qquad -D^{-2}e^{\mp_2it'D}\partial_1\partial_3a_1D^{-2}\partial_\nu e^{\mp_3it'D}\partial_1^2a_2)\,dt'[T^1,T^2]. \nonumber
\end{align}
Thus it suffices to consider the lower bound of $\mathcal F_x[\Lambda^1_\nu+\Lambda^2_\nu]$.
We consider $\nu=1$. First,
\begin{align*}
\mathcal F_x[\Lambda^1_1](t,\xi) & = \frac{1}{4i}\sum_{\pm_1,\pm_2,\pm_3}e^{\mp_1it|\xi|}\int \mathbf m_{123}(t,\xi,\eta)\frac{1}{|\xi||\xi-\eta|^2|\eta|^2} \\
& \qquad\qquad \times \left[(\xi_1-\eta_1)^2\eta_1^2\eta_2\chi_{-W_\lambda'}(\xi-\eta)\chi_{2W_\lambda}(\eta)-(\xi_1-\eta_1)(\xi_2-\eta_2)\eta_1^3\chi_{2W_\lambda}(\xi-\eta)\chi_{-W_\lambda'}(\eta)\right]\,d\eta,
\end{align*} 
where
$$
\mathbf m_{123}(t,\xi,\eta) = \frac{e^{it\omega_{123}}-1}{i\omega_{123}},\quad \omega_{123}=\pm_1|\xi|-\pm_2|\xi-\eta|-\pm_3|\eta|.
$$
The resonance case $\mathcal R$ occurs only when $\pm_1=\pm_2=\pm_3$ for which $|\omega_{123}|\lesssim\lambda^\frac12$. For other cases we have $|\omega_{123}|\approx\lambda$. We put $t=\epsilon\lambda^{-\frac12}$, for some $0<\epsilon\ll1$.
We also have
\begin{align*}
\mathcal F_x[\Lambda^2_1](t,\xi) & = \frac{1}{4i}\sum_{\pm_1,\pm_2,\pm_3}e^{\mp_1it|\xi|}\int \mathbf m_{123}(t,\xi,\eta)\frac{1}{|\xi||\xi-\eta|^2|\eta|^2} \\
& \qquad\qquad \times \left[(\xi_1-\eta_1)^2\eta_1^2\eta_3\chi_{-W_\lambda'}(\xi-\eta)\chi_{2W_\lambda}(\eta)-(\xi_1-\eta_1)(\xi_3-\eta_3)\eta_1^3\chi_{2W_\lambda}(\xi-\eta)\chi_{-W_\lambda'}(\eta)\right]\,d\eta.
\end{align*}
For $\eta\in-W_\lambda'$, $\eta_1$ is negative and hence $[\cdots]$ in the integrand is nonnegative.
Now we write
\begin{align*}
\mathcal F_x[\Lambda^1_1](t,\xi) 
&\ge \frac{1}{4}\left|\sum_{\pm_j\in\mathcal R}e^{\mp_1it|\xi|}\int\lambda^{-1}\lambda^\frac12t(1+O_{123}(\epsilon))\chi_{W_\lambda}(\xi-\eta)\chi_{-W_\lambda'}(\eta)\,d\eta\right| \\
& \qquad - \frac14\left|\sum_{\pm_j\in\mathcal R^c}\int\lambda^{-1}\lambda^{-1}\lambda^\frac12 \chi_{W_\lambda}(\xi-\eta)\chi_{-W_\lambda'}(\eta)\,d\eta\right| \\
&:= I + I\!I.
\end{align*}
Here, for $\xi-\eta\in{2W_\lambda},\eta\in{-W_\lambda'}$, we have $\xi\in{W_\lambda}$, and then for some $0<\rho<1$ we have $(1-\rho)\lambda<|\xi|<(1+\rho)\lambda$.
Clearly, we have $|I\!I|\lesssim \lambda^\frac12$. For $I$, we see that
\begin{align*}
I & \ge \frac14|\cos (t|\xi|)|\int t\lambda^{-\frac12}\chi_{W_\lambda}(\xi-\eta)\chi_{-W_\lambda'}(\eta)\,d\eta -\frac12 \left|\int tO_{123}(\epsilon)\lambda^{-\frac12}\chi_{W_\lambda}(\xi-\eta)\chi_{-W_\lambda'}(\eta)\,d\eta\right| \\
& := I_1 + I_2.
\end{align*}
Then we have $|I_2|\lesssim \epsilon^2\lambda$. 
Take $\lambda=\lambda(\epsilon,\rho,k)$ such that for each integer $k$,
$$
\frac{2k\pi-\epsilon}{\epsilon(1-\rho)}<\lambda^\frac12<\frac{2k\pi+\epsilon}{\epsilon(1+\rho)}.
$$
Then we have $|\cos(t|\xi|)|\approx1$ and hence $|I_1|\gtrsim \epsilon\lambda$. Thus we get $|\mathcal F_x[\Lambda^1_1](t,\xi)| \gtrsim \epsilon\lambda$ for $\xi\in W_\lambda$.
The treatment of $\Lambda^2_1$ is very similar. Since the integrand of $\Lambda^1_1$ and $\Lambda^2_1$ is both nonnegative, we obtain
$$
|\mathcal F_x[\Lambda^1_1+\Lambda^2_1](t,\xi)| \gtrsim \epsilon\lambda.
$$
On the other hand,
$$
\|\mathbf F\|_{H^r} \lesssim \|\partial_1a_2\|_{H^r}+\|\partial_2a_1\|_{H^r}+\|a_1a_2\|_{H^r} \lesssim 0 + \lambda^{r+\frac32}+\lambda^{r+1}.
$$
Thus for $C^2$ smoothness to hold, we must have
$
\epsilon\lambda^{s+2} \lesssim \lambda^{2r+2}+\lambda^{2r+3},
$
and hence
$
s-1 \le 2r.
$
\section*{Acknowledgements}
This research is supported by in part NRF-2018R1D1A3B07047782, NRF-2018R1A2B2006298, and NRF-2016K2A9A2A13003815. I would like to thank my mentor, Sebastian Herr, for bringing this problem to my attention. I am always grateful to all members of the IRTG 2235 in Bielefeld for their kind hospitality.


\end{document}